\documentclass[thmsa,11pt]{article}
\usepackage{latexsym}
\usepackage{amsmath,amssymb}

\newtheorem{theorem}{Theorem}[section]

\newtheorem{lemma}{Lemma}[section]
\newtheorem{proposition}{Proposition}[section]
\newtheorem{definition}{Definition}[section]

\newenvironment{proof}{\trivlist\item[\hskip
       \labelsep{\it Proof:\/}]\ignorespaces}{\hfill$\Box$\endtrivlist}

\input{tcilatex}

\begin{document}

\begin{center}
\textbf{A GODEL MODAL LOGIC \footnote{%
The results of this paper were announced at the meeting on "Logic,
Computability and Randomness", Cordoba, Argentina, Sept. 2004. Publication
was delayed, aiming to axiomatize the full logic with both modal operators,
which resulted elusive. Since the results have been quoted in some
publications based in an incomplete preliminary manuscript, we have chosen
to circulate this revision. We obtained recently the strong completeness of
the full logic,\ result which will appear elsewhere.}}

(revised Dec. 2008)\medskip

Xavier Caicedo $^{\ast }$

Ricardo Oscar Rodr\'{\i}guez $^{\ast \ast }$

\medskip

$^{\ast }${\small \ Departamento de Matem\'{a}ticas, Universidad de los
Andes, Bogot\'{a}, Colombia}

{\small xcaicedo@uniandes.edu.co }

$^{\ast \ast }${\small \ Departamento de Computaci\'{o}n, Fac. Ciencias
Exactas y Naturales }

{\small Universidad de Buenos Aires, 1428 Buenos Aires, Argentina}

{\small ricardo@dc.uba.ar }
\end{center}

\section{Introduction}

Sometimes it is needed in approximate reasoning to deal simultaneously with
both fuzziness of propositions and modalities, for instance one may try to
assign a degree of truth to propositions like \textquotedblleft John is
possibly tall\textquotedblright\ or \textquotedblleft John is necessarily
tall\textquotedblright , where \textquotedblleft John is
tall\textquotedblright\ is presented as a fuzzy proposition. Fuzzy logic
should be a suitable tool to model not only vagueness but also other kinds
of information features like certainty, belief or similarity, which have a
natural interpretation in terms of modalities.

We address in this paper the case of pure modal operators for G\"{o}del
logic, one of the main systems of fuzzy logic arising from H\'{a}jek \cite%
{Hajek98} classification. For this purpose we consider a many-valued version
of Kripke semantics for modal logic where both, propositions at each world
and the accessibility relation, are infinitely valued in the standard G\"{o}%
del algebra [0,1].

We provide strongly complete axiomatizations for the $\square $-fragment and
the $\Diamond $-fragment of the resulting minimal logic. These fragments are
shown to behave quite asymmetrically. Validity\ in the first one is
univocally determined by the class of frames having a crisp (that is,
two-valued)\ accessibility relation, while validity in the second requires
truly fuzzy frames. In addition, the $\square $-fragment does not enjoy the
finite model property with respect to the number of worlds or the number of
truth values while the $\Diamond $-fragment does.

We consider also the G\"{o}del analogues of the classical modal systems T,
S4 and S5 for each modal operator and show that the first two are
characterized by the many-valued versions of the frame properties\ which
characterize their classical counterparts.

Our approach is related to Fitting \cite{Fitting92} who considers Kripke
models taking values in a fixed finite Heyting algebra;\ however, his
systems and completeness proofs depend essentially on finiteness of the
algebra and he fact that his languages contain constants for all the truth
values of the algebra. We most relay on completely different methods.

Modal logics with an intuitionistic basis and Kripke style semantics have
been investigated in a number of relevant papers (see\ Ono \cite{Ono77},
Fischer Servi \cite{FischerServi84}, Bo\u{s}zic and Do\u{s}en \cite{Dosen84}%
, Font \cite{Font86}, Wolter \cite{Wolter97}, from an extensive literature),
but in all cases the models carry two (or more)\ crisp accessibility
relations satisfying some commuting properties: a pre-order to account for
the intuitionistic connectives\ and one or more binary relations to account
for the modal operators. Our semantics has,.instead, a single arbitrary
fuzzy accessibility relation and does not seem reducible to those
multi-relational semantics since the latter enjoy the finite model property
for $\square $ (cf.Grefe \cite{Grefe98}).

We assume the reader is acquainted with modal and G\"{o}del logics and the
basic laws of linear Heyting algebras (cf. Chagrov \cite{Chagrov97}).

\section{G\"{o}del-Kripke models}

The language\textit{\ }$\mathcal{L}_{\square \Diamond }$ of propositional 
\emph{G\"{o}del modal logic} is built from a set $Var$ of propositional
variables, logical connectives symbols $\wedge ,\rightarrow ,\bot ,$ and the
modal operator symbols $\square $ and $\Diamond $. Other connectives are
defined:\medskip

$\top :=\varphi \rightarrow \varphi $

$\lnot \varphi :=\varphi \rightarrow \bot $

$\varphi \vee \psi :=((\varphi \rightarrow \psi )\rightarrow \psi ))\wedge
((\psi \rightarrow \varphi )\rightarrow \varphi ))$

$\varphi \longleftrightarrow \psi :=(\varphi \rightarrow \psi )\wedge (\psi
\rightarrow \varphi ).$\medskip

\noindent $\mathcal{L}_{\square }$ and $\mathcal{L}_{\Diamond }$ will
denote, respectively, the $\square $-fragment and the $\Diamond $-fragment
of the language.

As stated before, the semantics of G\"{o}del modal logic will be based in
fuzzy Kripke models where the valuations at each world and also the
accessibility relation between worlds are $[0,1]$-valued. The symbols $\cdot 
$ and $\Rightarrow $ will denote the G\"{o}del norm in $[0,1]$ and its
residuum, respectively: 
\begin{equation*}
\begin{array}{ll}
a\cdot b=\min \{a,b\},\text{ \ \ \ \ \ } & a\Rightarrow b=\left\{ 
\begin{array}{ll}
1, & \mathrm{if}\ a\leq b \\ 
b, & \mathrm{otherwise}%
\end{array}%
\right. 
\end{array}%
\end{equation*}%
the derived maximum and pseudo-complement operations will be denoted $%
\curlyvee $ and $-,$ respectively. This yields the standard G\"{o}del
algebra; that is, the unique Heyting algebra structure in the linearly
ordered interval.

\begin{definition}
\label{models} A G\"{o}del-Kripke model (GK-model)\ \emph{will be a
structure }$\langle W,S,e\rangle $\emph{\ where:}
\end{definition}

\noindent $\bullet $ $\ W$\ is a non-empty set of objects that we call \emph{%
worlds} of $M.$

\noindent $\bullet $ $\ S:W\times W\rightarrow \lbrack 0,1]$ is an
arbitrary\ function $(x,y)\longmapsto Sxy$.

\noindent $\bullet $ $\ e:W\times Var\rightarrow \lbrack 0,1]$ is an
arbitrary\ function $(x,p)\longmapsto e(x,p)$.\medskip

The evaluations $e(x,-):Var\rightarrow \lbrack 0,1]$\ are extended
simultaneously to all formula in $\mathcal{L}_{\square \Diamond }$\ by
defining inductively at each world $x$:$\ $

\medskip

$e(x,\varphi \wedge \psi ):=e(x,\varphi )\cdot e(x,\psi )$

$e(x,\varphi \rightarrow \psi ):=e(x,\varphi )\Rightarrow e(x,\psi )$

$e(x,\bot ):=0$

$e(x,\Box \varphi ):=\inf_{y\in W}\{Sxy\Rightarrow e(y,\varphi )\}$

$e(x,\Diamond \varphi ):=\sup_{y\in W}\{Sxy\cdot e(y,\varphi )\}$.

\medskip

\noindent It follows that $e(x,\varphi \vee \psi )=e(x,\varphi )\curlyvee
e(x,\psi )$ and $e(x,\lnot \varphi )=-e(x,\varphi ).$\medskip

The notions of a formula $\varphi $\ being true at a world $x$, valid in a
model $M=\langle W,S,e\rangle ,$ or universally valid, are the usual ones:%
\emph{\medskip }

$\varphi $ is \emph{true in }$M$ \emph{at} $x$, written $M\models
_{x}\varphi ,$ iff $e(x,\varphi )=1$.

$\varphi $ is \emph{valid in} $M$, written $M\models \varphi ,$ iff $%
M\models _{x}\varphi $ at any world $x$ of $M.$

$\varphi $ is GK-\emph{valid}, written $\models _{GK}\varphi $, if it is
valid in all the GK-models.\emph{\medskip }

\noindent Clearly, all valid schemes of G\"{o}del logic are GK-valid. In
addition,

\begin{proposition}
\label{soundness}. The following modal schemes are\textbf{\ }GK-valid:
\end{proposition}

\noindent $%
\begin{array}{lll}
\mathbf{K}_{\square } &  & \square (\varphi \rightarrow \psi )\rightarrow
(\square \varphi \rightarrow \square \psi ) \\ 
\mathbf{Z}_{\square } &  & \lnot \lnot \square \theta \rightarrow \square
\lnot \lnot \theta \\ 
\mathbf{D}_{\Diamond } &  & \Diamond (\varphi \vee \psi )\rightarrow
(\Diamond \varphi \vee \Diamond \psi )\text{ \ (in fact, an equivalence)} \\ 
\mathbf{Z}_{\Diamond } &  & \Diamond \lnot \lnot \varphi \rightarrow \lnot
\lnot \Diamond \varphi \\ 
\mathbf{F}_{\Diamond } &  & \lnot \Diamond \bot%
\end{array}%
$

\proof Let $M=\langle W,S,e\rangle $ be an arbitrary G\"{o}del-Kripke model
and $x\in W.$\emph{\medskip }

\noindent ($\mathbf{K}_{\square })$ By definition \ref{models} and
properties of the residuum we have for any $y\in W$: $e(x,\square (\varphi
\rightarrow \psi ))\cdot e(x,\square \varphi )$ $\leq (Sxy\Rightarrow
(e(y,\varphi )\Rightarrow e(y,\psi ))\cdot (Sxy\Rightarrow e(y,\varphi ))$ $%
\leq (Sxy\Rightarrow e(y,\psi )).$ Taking the meet over $y$ in the last
expression: $e(x,\square (\varphi \rightarrow \psi ))\cdot e(x,\square
\varphi )$ $\leq e(x,\square \psi ),$ hence\ $e(x,\square (\varphi
\rightarrow \psi ))$ $\leq e(x,\square \varphi \rightarrow \square \psi )$.%
\emph{\medskip }

\noindent ($\mathbf{Z}_{\square })$ Utilizing the Heyting algebra identity: $%
--(x\Rightarrow y)=(x\Rightarrow --y),$ we have: $e(x,\lnot \lnot \Box
\theta )=--e(x,\Box \theta )$ $\leq --(Sxy\Rightarrow e(y,\theta ))$ $%
=(Sxy\Rightarrow --e(y,\theta ))=(Sxy\Rightarrow e(y,\lnot \lnot \theta )).$
Taking meet over $y$ in the last expression: $e(x,\lnot \lnot \Box \theta
)\leq e(x,\Box \lnot \lnot \theta ).$

\noindent ($\mathbf{D}_{\Diamond }$)\ By properties of suprema and
distributivity of $\curlyvee $ over $\cdot ,$ $e(\Diamond (\varphi \vee \psi
))$ $=\sup_{y}\{Sxy\cdot (e(y,\varphi )\curlyvee e(y,\varphi ))\}$ $%
=\sup_{y}\{Sxy\cdot e(y,\varphi )\}\curlyvee \sup_{y}\{Sxy\cdot e(y,\varphi
)\}$

\noindent ($\mathbf{Z}_{\Diamond })$ $Sxy\cdot e(\lnot \lnot \varphi ,y)\leq
--(Sxy\cdot e(\varphi ,y))\leq --e(\Diamond \varphi ,x)=e(\lnot \lnot
\Diamond \varphi ,x).$\emph{\medskip }

\noindent ($\mathbf{F}_{\Diamond })$ $e(x,\Diamond \bot )=\sup_{y}\{Sxy\cdot
0\}=0.$ \ $\blacksquare $\emph{\medskip }

The Modus Ponens rule preserves truth at every world of any GK-model. On the
other hand, the classical introduction rules for the modal operators 
\begin{equation*}
\begin{array}{lllll}
\mathbf{RN}_{\square }:%
\begin{array}{l}
\underline{\varphi } \\ 
\square \varphi%
\end{array}
&  &  &  & \mathbf{RN}_{\Diamond }:%
\begin{array}{l}
\underline{\varphi \rightarrow \psi } \\ 
\Diamond \varphi \rightarrow \Diamond \psi%
\end{array}%
.%
\end{array}%
\end{equation*}%
do not preserve local truth. However,

\begin{proposition}
\label{soundness2}$\mathbf{RN}_{\square }$ and $\mathbf{RN}_{\Diamond }$
preserve validity at any given model, thus they preserve GK-validity.
\end{proposition}

\proof. ($\mathbf{RN}_{\square })$ If $e(x,\varphi )=1$ for all $x$ then $%
e(x,\square \varphi )=\inf_{y}\{Sxy\Rightarrow e(y,\varphi )\}=\inf \{1\}=1$
for all $x$. ($\mathbf{RN}_{\Diamond })$ If $e(x,\varphi \rightarrow \psi
)=1 $ for all $x$ then $Sxy\cdot e(y,\varphi )\leq Sxy\cdot e(y,\psi )\leq
e(x,\Diamond \psi ).$ Taking join over $y$ in the left hand side of the last
inequality, $e(x,\Diamond \varphi )\leq e(x,\Diamond \psi ).$ $\
\blacksquare $\emph{\medskip }

Semantic consequence is defined for any theory $T\subseteq \mathcal{L}%
_{\square \Diamond },\ $as follows:

\begin{definition}
\label{consequence}$T\models _{GK}\varphi $ if and only if for any GK-model $%
\mathcal{M}$ and any world $x$ in $M,$ $\mathcal{M}\models _{x}T$ implies $%
\mathcal{M}\models _{x}\varphi .$
\end{definition}

An alternative notion of logical consequence arises naturally. Set $%
e(x,T)=\{e(x,\varphi ):\varphi \in T\}$ then:

\begin{definition}
\label{ConseqInf}$T\models _{GK\leq }\varphi $\textit{\ if and only if for
any GK-model }$M$\textit{\ and any world }$x$\textit{\ in }$M,$\textit{\ }$%
\inf e(x,T)\leq e(x,\varphi ).$
\end{definition}

Clearly, $\models _{GK\leq }$ implies $\models _{GK}$, and it will follow
from our completeness theorems that both notions are equivalent for
countable theories. This fact has been already observed for pure G\"{o}del
logic by Baaz and Zach in \cite{BaazZach98}.

Note that Modus Ponens preserves consequence but this is not the case of the
inference rules $\mathbf{RN}_{\square }$ and $\mathbf{RN}_{\Diamond }.$

\section{On strong completeness of G\"{o}del logic}

\label{GodelCompl}

To prove strong completeness of the unimodal fragments $\mathcal{L}_{\square
}$ and $\mathcal{L}_{\Diamond }$ we will reduce the problem to pure G\"{o}%
del propositional logic.

In the rest of this paper $\mathcal{L}(X)$ will denote the G\"{o}del
language built from a set of propositional variables $X$ and the connectives 
$\wedge ,\rightarrow ,\bot .$

Let $\mathcal{G}$ be a fixed axiomatic calculus for G\"{o}del logic, say the
following one given by H\'{a}jek (\cite{Hajek98}, Def. 4.2.3.):

\medskip

$(\varphi \rightarrow \psi )\rightarrow ((\psi \rightarrow \chi )\rightarrow
(\varphi \rightarrow \chi ))$

$(\varphi \wedge \psi )\rightarrow \varphi $

$(\varphi \wedge \psi )\rightarrow (\psi \wedge \varphi )$

$(\varphi \rightarrow (\psi \rightarrow \chi ))\longleftrightarrow ((\varphi
\wedge \psi )\rightarrow \chi )$

$((\varphi \wedge \psi )\rightarrow \chi )\longleftrightarrow (\varphi
\rightarrow (\psi \rightarrow \chi ))$

$\varphi \rightarrow (\varphi \wedge \varphi )$

$((\varphi \rightarrow \psi )\rightarrow \chi )\rightarrow (((\psi
\rightarrow \varphi )\rightarrow \chi )\rightarrow \chi )$

$\bot \rightarrow \varphi $\medskip

\textbf{MP}: From $\varphi $ and $\varphi \rightarrow \psi $, infer $\psi $

\medskip

$\vdash $ will denote deduction in this calculus.\medskip

\noindent It is well known that $\mathcal{G}$ is deductively equivalent to
Dummett logic, the intermediate logic obtained by adding to Heyting calculus
the pre-linearity schema: \medskip

$(\varphi \rightarrow \psi )\vee (\psi \rightarrow \varphi )$. \medskip

Given a valuation $v:X\rightarrow \lbrack 0,1],$ let $\overline{v}$ denote
the extension of $v$ to $\mathcal{L}(X)$ according to the G\"{o}del
interpretation of the connectives. We will need the following strong form of
standard completeness for G\"{o}del logic:

\begin{proposition}
{\label{Godel}}Let $T$\ be a countable\ theory and $U$\ a countable set of
formulas of $\mathcal{L}(X)$ such that for every finite $S\subseteq U$\ we
have $T\nvdash \bigvee S$\ then there is a valuation $v:X\rightarrow \lbrack
0,1]$\ such that $\overline{v}(\alpha )=1$\ for all $\alpha \in T$\ and $%
\overline{v}(\beta )<1$ for each $\beta \in U.$
\end{proposition}

\proof Extend $T$ to a prime theory $T^{\prime }$ (that is, $T^{\prime
}\vdash \alpha \vee \beta $ implies $T^{\prime }\vdash \alpha $ or $%
T^{\prime }\vdash \beta $ ) satisfying the same hypothesis with respect to $%
U $ (this is standard)$.$ The Lindenbaum algebra $\mathcal{L}(X)/_{\equiv
T^{\prime }}$ of $T^{\prime }$ is linearly ordered since by primality and
the pre-linearity schema $T^{\prime }\vdash \alpha \rightarrow \beta $ or $%
T^{\prime }\vdash \beta \rightarrow \alpha $. Moreover, the valuation $%
v:X\rightarrow \mathcal{L}(X)/_{T^{\prime }},$ $v(x)=x/_{\equiv T^{\prime }}$
is such that $v(T)=1,$ $v(\beta )<1$ for all $\beta \in U.$ As $T^{\prime }$
is countable we may assume $X$ is countable and thus, being also countable, $%
\mathcal{L}(X)/_{\equiv T^{\prime }}$ is embeddable in the G\"{o}del algebra 
$[0,1],$ therefore, we may assume $v:X\rightarrow \lbrack 0,1]$. $%
\blacksquare $

\medskip

From the proposition we obtain the usual formulation of completeness for
countable $T$. We can not expect strong standard completeness\ of\ $\mathcal{%
G}$ for uncountable theories, as the following example illustrates.

\medskip

\noindent \textbf{Example}. Set $T=\{(p_{\beta }\rightarrow p_{\alpha
})\rightarrow q:\alpha <\beta <\omega _{1}\}$ where $\omega _{1}$ is the
first uncountable cardinal, then $T\nvdash q.$ Otherwise we would have $%
\Sigma \vdash q,$ for some finite $\Sigma =\{(p_{\alpha _{i+1}}\rightarrow
p_{\alpha _{i}})\rightarrow q:1\leq i<n\},$ but this is not possible by
soundness of $\mathcal{G},$ because the valuation $v(q)=\frac{1}{2}$, $%
v(p_{\alpha _{i}})=\frac{1}{2}(1-\frac{1}{i+1})$ for $1\leq i\leq n,$ makes $%
v(p_{\alpha _{i}})<v(p_{\alpha _{i+1}})<\frac{1}{2}$ and thus $\overline{v}%
((p_{\alpha _{i+1}}\rightarrow p_{\alpha _{i}})\rightarrow q)=1$ for $1\leq
i<n,$ while $v(q)<1.$ On the other hand, there is no valuation $v$ such that 
$\overline{v}(T)=1$ and $v(q)<1$, because that would imply $\overline{v}%
(p_{\beta }\rightarrow p_{\alpha })<1$ for all $\alpha <\beta <\omega _{1},$
and thus the set $\{v(p_{\alpha }):\alpha <\omega _{1}\}$ would be ordered
in type $\omega _{1},$ which is impossible because any well ordered subset
of $([0,1],<)$ is at most countable.

\section{Completeness of the $\square $-fragment}

Let $\mathcal{G}_{\square }$ be the formal system on the language $\mathcal{L%
}_{\square }$ which is obtained by adding to the system $\mathcal{G}$ for G%
\"{o}del logic (applied to $\mathcal{L}_{\square }$) the following axiom
schemes and rule: \medskip

$\mathbf{K}_{\square }$: \ $\square (\varphi \rightarrow \psi )\rightarrow
(\square \varphi \rightarrow \square \psi )$

$\mathbf{Z}_{\square }$: $\ \lnot \lnot \square \theta \rightarrow \square
\lnot \lnot \theta $

$\mathbf{NR}_{\square }$: From $\varphi $ infer $\square \varphi $\medskip

$\vdash _{\mathcal{G}_{\square }}\varphi $ expresses theorem-hood in this
logic. Proofs with assumptions are also allowed, with the restriction that $%
\mathbf{NR}_{\square }$ is to be applied to theorems only (or, what amounts
to the same, to previous steps of the proof not depending on the
assumptions). $T\vdash _{\mathcal{G}_{\square }}\varphi $ will express that
there is such a proof of $\varphi $ with assumptions from the set $T.$

The deduction theorem follows readily by induction in the length of
proofs:\medskip

$\mathbf{DT}$: $\ T\cup \{\alpha \}\vdash _{\mathcal{G}_{\square }}\varphi $%
\emph{\ implies} $T\vdash _{\mathcal{G}_{\square }}\alpha \rightarrow
\varphi .$

\medskip

\noindent Applying consecutively $\mathbf{DT}$, $\mathbf{NR}_{\square }$, $%
\mathbf{K}_{\square }$, and $MP,$ we obtain the derived rule:

\begin{lemma}
{\label{necessit}}If $\mu _{1},...,\mu _{k}\vdash _{\mathcal{G}_{\square
}}\varphi $ then $\Box \mu _{1},...,\Box \mu _{k}\vdash _{\mathcal{G}%
_{\square }}\Box \varphi .$
\end{lemma}

\noindent We obtain also soundness of $\mathcal{G}_{\square }$:

\begin{lemma}
{\label{sound}}$T\vdash _{\mathcal{G}_{\square }}\varphi $ \ implies $\
T\models _{GK\leq }\varphi $, hence, $T\models _{GK}\varphi .$
\end{lemma}

\proof By the deduction theorem, $T\vdash _{\mathcal{G}_{\square }}\varphi $
implies $\vdash _{\mathcal{G}_{\square }}(\wedge \Sigma \rightarrow \varphi
) $ for some finite $\Sigma \subseteq T.$ Since the axioms of $\mathcal{G}%
_{\square }$ are valid in all GK-models (Prop. \ref{soundness}) and $\mathbf{%
MP},$ $\mathbf{NR}_{\square }$, $\mathbf{K}_{\square }$ preserve validity
(Prop. \ref{soundness2}) then $\models _{GK}(\wedge \Sigma \rightarrow
\varphi ).$ Therefore, $\inf e(x,T)\leq e(x,\wedge \Sigma )\leq e(x,\varphi
) $ for any world $x$ in any GK-model. $\ \blacksquare $

\medskip

Let 
\begin{equation*}
T\mathcal{G}_{\square }=\{A:A\text{ is a theorem of }\mathcal{G}_{\square }\}%
\text{ }
\end{equation*}%
Since all uses of $\mathbf{NR}_{\square }$ in a proof of $T\vdash _{\mathcal{%
G}_{\square }}\varphi $ produce theorems of $\mathcal{G}_{\square }$, the
proof may bee seen as one in which Modus Ponens is the only rule utilized
and $T\mathcal{G}_{\square }$ is part of the assumptions. That is,

\begin{lemma}
{\label{A}}. $T\vdash _{\mathcal{G}_{\square }}\varphi $ if and only if $%
T\cup T\mathcal{G}_{\square }\vdash \varphi $ in pure G\"{o}del logic.
\end{lemma}

To prove strong completeness of $\mathcal{G}_{\square }$ we will define a
canonical GK-model with the property that for any countable theory $T$ and
any formula $\varphi \,$such that $T\nvdash _{\mathcal{G}_{\square }}\varphi
,$ there is a world $x$ in the model which assigns the value 1 to $T$ but
less than 1 to $\varphi $. A surprising fact will be that this may be
achieved with a model where the accessibility relation is crisp.

Let $\square \mathcal{L}_{\square }=\{\square \theta :\theta \in \mathcal{L}%
_{\square }\}$ be the set of formulas in $\mathcal{L}_{\square }$ which
start with the connective $\square $. Then any formula in $\mathcal{L}%
_{\square }$ may be seen as a formula of the pure G\"{o}del language built
from $X=Var\cup \square \mathcal{L}_{\square }$ by means of $\wedge ,\lnot
,\bot .$ That is, we may consider the formulas in $\square \mathcal{L}%
_{\square }$ as additional propositional variables for G\"{o}del
logic.\medskip

\textbf{Canonical model }$\mathcal{M}_{\square }=(W^{\ast },S^{\ast
},e^{\ast })$.\medskip

\noindent $\bullet $ $\ $The set of worlds $W^{\ast }$ will consist of those
valuations $v:Var\cup \square \mathcal{L}_{\square }\rightarrow \lbrack 0,1]$
which satisfy $\overline{v}(T\mathcal{G}_{\square })=1$ when extended to $%
\overline{v}:\mathcal{L}_{\square }=\mathcal{L}(Var\cup \square \mathcal{L}%
_{\square })\rightarrow \lbrack 0,1]$ according to the G\"{o}del
interpretation of $\wedge ,\rightarrow ,\bot .$

\noindent $\bullet $ $\ $The fuzzy accessibility relation between worlds in $%
\mathcal{M}_{\square }$ (actually a crisp relations)\ will be given by 
\begin{equation*}
S^{\ast }vw=\left\{ 
\begin{array}{ll}
1, & \text{ if }v(\Box \theta )\leq \overline{w}(\theta ),\text{ for all }%
\theta \in \mathcal{L}_{\square } \\ 
0, & \mathrm{otherwise}%
\end{array}%
\right. ,
\end{equation*}%
\noindent $\bullet $ $\ $The valuation associated to the world $v$ will be $%
v\upharpoonright Var.$ That is, $e^{\ast }(p,v)=v(p)$ for any $p\in Var.$%
\medskip

For the sake of simplicity, we will write $v(\varphi )$ for $\overline{v}%
(\varphi ),$ from now on.\medskip

\begin{lemma}
{\label{equation}}\emph{For any world }$v$ \emph{in the canonical model }$%
\mathcal{M}_{\square }$\emph{\ and any} $\varphi $ 
\begin{equation*}
e^{\ast }(v,\varphi )=v(\varphi ).
\end{equation*}
\end{lemma}

\proof This is proven by induction in the complexity of $\varphi $ seen
again as a formula of $\mathcal{L}_{\square }.$ The atomic step and the
inductive steps for the G\"{o}del connectives being straightforward, it is
enough to verify inductively $e^{\ast }(v,\Box \varphi )=v(\Box \varphi )$.
By induction hypothesis we may assume $e^{\ast }(w,\varphi )=w(\varphi )$
for any $w,$and thus we must show%
\begin{equation*}
v(\Box \varphi )=\inf_{w}\{w(\varphi ):S^{\ast }vw=1\}
\end{equation*}%
By definition, $S^{\ast }vw=1$ implies $v(\Box \varphi )\leq w(\varphi ),$
hence 
\begin{equation*}
v(\Box \varphi )\leq \inf_{w}\{w(\varphi ):S^{\ast }vw=1\}.
\end{equation*}%
Since equality above is trivial for $v(\square \varphi )=1$, it remains only
to show in case $v(\square \varphi )=\alpha <1$ that%
\begin{equation}
\inf_{w}\{w(\varphi ):S^{\ast }vw=1\}\leq \alpha .  \label{equalnec1}
\end{equation}%
That is, for any $\epsilon >0$ there is $w$ such that $S^{\ast }vw=1$ and $%
w(\varphi )<\alpha +\epsilon $. To see this we prove first:\medskip\ 

\noindent \textbf{Claim}{. }\emph{Let }$v$ \emph{be a world of }$\mathcal{M}%
_{\square }$\emph{\ and }$\varphi $\emph{\ be such that }$v(\Box \varphi
)=\alpha <1,$\emph{\ then there exists a world }$u$\emph{\ of }$\mathcal{M}%
_{\square }$\emph{\ such that }$u(\varphi )<1$ \emph{and}$\newline
($i) \ $u(\theta )=1$\emph{\ if }$v(\Box \theta )>\alpha \newline
($ii) $u(\theta )>0$\emph{\ if }$v(\Box \theta )>0.$\medskip

\proof Assume $v(\Box \varphi )=\alpha <1\ $and set 
\begin{equation*}
T_{\varphi ,v}=\{\theta :v(\Box \theta )>\alpha \}\cup \{\lnot \lnot \theta
:v(\Box \theta )>0\}
\end{equation*}%
Notice that $v(\Box \mu )>\alpha $ for any $\mu \in T_{\varphi ,v}$ because $%
v(\Box \theta )>0$ implies $v(\lnot \lnot \Box \theta )=1,$ and thus $v(\Box
\lnot \lnot \theta )=1$ since $v$ satisfies axiom $\mathbf{Z}_{\square }.$
This implies that $T_{\varphi ,v}\nvdash _{\mathcal{G}_{\square }}\varphi .$
Otherwise, $\mu _{1},...,\mu _{k}\vdash _{\mathcal{G}_{\square }}\varphi $
for some $\mu _{i}\in T_{\varphi ,v}$ and thus 
\begin{equation*}
\Box \mu _{1},...,\Box \mu _{k}\vdash _{\mathcal{G}_{\square }}\Box \varphi
\end{equation*}%
by Lemma \ref{necessit}. Hence, by Lemma \ref{sound} and the previous
observations, 
\begin{equation*}
\alpha <\min \{\Box \mu _{1},...,\Box \mu _{k}\}\leq v(\Box \varphi ),
\end{equation*}%
a contradiction. By Lemma \ref{A} we have $T_{\varphi ,v}\cup T\mathcal{G}%
_{\square }\nvdash _{{}}\varphi $ and by countability of $T_{\varphi ,v}\cup
T\mathcal{G}_{\square }$ we may use the completeness theorem of G\"{o}del
logic (Proposition \ref{Godel}) to get a G\"{o}del valuation $u:L\rightarrow
\lbrack 0,1]$ such that $u(T_{\varphi ,v})=1$ and $u(\varphi )<1.$ Then $%
u\in \mathcal{M}_{\square }$ and (i) holds by construction. Moreover, (ii)
is satisfied because $u(\lnot \lnot \theta )=1$ and thus $u(\theta )>0$ if%
\emph{\ }$v(\Box \theta )>0.$ This ends the proof of the claim.\medskip

Pick now an strictly increasing function $g:[0,1]\rightarrow \lbrack 0,1]$
such that 
\begin{equation*}
g(1)=1,\,\text{\ }g(0)=0,\,\text{\ and }g[(0,1)]=(\alpha ,\alpha +\epsilon ).
\end{equation*}%
As $g$ is an homomorphism of Heyting algebras, the valuation $w=g\circ u$
preserves the value 1 of the formulas in $T\mathcal{G}_{\square }$ and thus
it belongs to $\mathcal{M}_{\square }$. Moreover, $v(\Box \theta )\leq
w(\theta )$ for all $\theta $:\ 

- if $v(\Box \theta )>\alpha $ because $w(\theta )=g(u(\theta ))=g(1)=1\ $by
(1) above.

- if $0<v(\Box \theta )\leq \alpha $ because then $0<u(\theta )\leq 1$ by
(2) above, and thus $w(\theta )$ $=g(u(\theta ))$ $\in (\alpha ,\alpha
+\epsilon )\cup \{1\}.$

\noindent This means $S^{\ast }vw=1$, and since $u(\varphi )<1$ we have, $%
w(\varphi )=g(u(\varphi ))<\alpha +\epsilon ,$ which shows \ref{equalnec1}. $%
\blacksquare $

\bigskip

Call a GK-model \emph{crisp} if $S:W\times W\rightarrow \{0,1\},$ and write $%
T\models _{Crisp}\varphi $ if the consequence relation holds at each node of
any crisp GK-model.

\begin{theorem}
{\label{completnec}}For any countable theory $T$ and formula $\varphi $ in $%
\mathcal{L}_{\square }$ the following are equivalent:
\end{theorem}

\noindent $%
\begin{array}{lll}
\text{(i)} &  & T\vdash _{\mathcal{G}_{\square }}\varphi \\ 
\text{(ii)} &  & T\models _{GK_{\leq }}\varphi \\ 
\text{(iii)} &  & T\models _{GK} \\ 
\text{(iv)} &  & T\models _{Crisp}\varphi .%
\end{array}%
$

\proof By Lemma \ref{sound}, it is enough to show (iv) $\Rightarrow $ (i).
If $T\nvdash _{\mathcal{G}_{\square }}\varphi $ then $T\cup T\mathcal{G}%
_{\square }\nvdash \varphi $ by Lemma \ref{A}, and by strong completeness of
G\"{o}del logic there is a valuation $v:Var\cup \square \mathcal{L}_{\square
}\rightarrow \lbrack 0,1]$ such that $\overline{v}(T)=\overline{v}(T\mathcal{%
G}_{\square })=1$ and $\overline{v}(\varphi )<1.$ Hence, $v\in W^{\ast }$ by
definition, $e^{\ast }(v,T)=\overline{v}(T)=1,$ and $e^{\ast }(v,\varphi )=%
\overline{v}(\varphi )<1$ by Lemma \ref{equation}, showing that $\mathcal{M}%
_{\square }\models _{v}T$ but $\mathcal{M}_{\square }\nvDash _{v}\varphi .$
That is, $T\nvDash _{Crisp}\varphi $ because the canonical model is crisp. $%
\blacksquare $

\medskip

The example in Section \ref{GodelCompl} shows that standard strong
completeness does not hold in modal G\"{o}del logic with respect to
uncountable theories.

\section{$\mathcal{G}_{\square }$ does not have the finite model property}

\label{finitemodel}

The following example shows that $\mathcal{G}_{\square }$ does not have the
finite model property with respect to GK-models. The scheme%
\begin{equation*}
\square \lnot \lnot \theta \rightarrow \lnot \lnot \square \theta ,
\end{equation*}%
reciprocal of axiom $\mathbf{Z}_{\square }$, is not valid because it fails
in the (crisp) model $\mathcal{M}$ = $(\mathbb{N},S,e),$ where%
\begin{eqnarray*}
Smn &=&1\text{ for all }m,n \\
e(n,p) &=&\tfrac{1}{n+1}\text{ for all }n
\end{eqnarray*}%
Indeed, $e(n,\lnot \lnot p)=--\frac{1}{n+1}=1$ for all $n$ and thus, $%
e(0,\Box \lnot \lnot p)=\inf \{1\Rightarrow 1\}=1.$ On the other hand, $%
e(0,\Box p)=\inf_{n\in \mathbb{N}}\{1\Rightarrow \frac{1}{n+1}\}=0,$ and
thus $e(0,\lnot \lnot \Box p)=0.$

However,

\begin{theorem}
$\Box \lnot \lnot \theta \rightarrow \lnot \lnot \Box \theta $ is valid in
any GK-model $\langle W,S,e\rangle $\ with\ finite $W$.
\end{theorem}

\proof Given a model $\mathcal{M}$ = $\langle W,S,e\rangle ,$ we have:%
\begin{equation}
e(v,\Box \lnot \lnot \theta )=\left\{ 
\begin{array}{ll}
0, & \exists w\in W:Svw>0\ \text{\ and \ }e(w,\theta )=0 \\ 
1, & \text{otherwise}%
\end{array}%
\right.  \label{menos}
\end{equation}%
Now, $e(v,\lnot \lnot \Box \theta )=0$ iff and only if $e(v,\Box \theta )=0,$
which means there is a sequence of worlds $\{w_{n}\}_{n}$ such that $%
\{Svw_{n}\Rightarrow e(w_{n},\theta )\}_{n}$ converges to $0,$ that is,%
\begin{equation}
e(\lnot \lnot \Box \theta ,v)=\left\{ 
\begin{array}{ll}
0, & \text{if}%
\begin{array}{l}
\exists \{w_{n}\}\subseteq W:Svw_{n}>e(w_{n},\theta )\text{ for all } \\ 
n\in \mathbb{N}\text{, and }\{e(w_{n},\theta )\}_{n}\text{ converges to }0%
\end{array}
\\ 
1, & \text{otherwise}%
\end{array}%
\right.  \label{mas}
\end{equation}%
Assume $e(v,\Box \lnot \lnot \theta )=1.$ Then, according to (\ref{menos})$,$
$Svw=0$ or $e(w,\theta )>0$ for any $w\in W.$ If we had $e(v,\lnot \lnot
\Box \theta )=0,$ then the sequence $\{w_{n}\}_{n\in \mathbb{N}}$ given by (%
\ref{mas}) would satisfy:\ $Svw_{n}>e(w_{n},\theta ),$ and hence $%
e(w_{n},\theta )>0$ for all $n$ by the previous observation. If $W$ is
finite, the set $\{e(w_{n},\theta ):n\in \mathbb{N}\}$ has a minimum
positive value and thus the sequence, $\{e(w_{n},\theta )\}_{n}$ would not
converge to $0,$ a contradiction. $\blacksquare $

\medskip

The proof of the theorem shows that $\Box \lnot \lnot \theta \rightarrow
\lnot \lnot \Box \theta $ would be valid in all GK models with finite G\"{o}%
del algebra of values.

\section{Completeness of the $\Diamond $-fragment}

The system $\mathcal{G}_{\Diamond }$ results by adding to $\mathcal{G}$ the
following axiom schemes and rule in the language $\mathcal{L}_{\Diamond }:$%
\medskip

\textbf{D}$_{\Diamond }$: $\ \Diamond (\varphi \vee \psi )\rightarrow
(\Diamond \varphi \vee \Diamond \psi )$

\textbf{Z}$_{\Diamond }$: $\ \ \Diamond \lnot \lnot \varphi \rightarrow
\lnot \lnot \Diamond \varphi $

\textbf{F}$_{\Diamond }$: $\ \ \lnot \Diamond \bot $

\textbf{RN}$_{\Diamond }$: \ \textit{From }$\varphi \rightarrow \psi $ 
\textit{infer} $\Diamond \varphi \rightarrow \Diamond \psi $\medskip

As in the case of the $\square $-fragment, in proofs with assumptions the
rule \textbf{RN}$_{\Diamond }$ is to be used in theorems only, and under
this definition we have the deduction theorem $\mathbf{DT}$, the derived
rule:

\begin{lemma}
{\label{possibilit}}If $\varphi \vdash _{\mathcal{G}_{\Diamond }}\psi $ 
\textit{then} $\Diamond \varphi \vdash _{\mathcal{G}_{\Diamond }}\Diamond
\psi .$
\end{lemma}

\noindent and the soundness theorem:

\begin{lemma}
{\label{sound2}}$T\vdash _{\mathcal{G}\Diamond }\varphi $ \ implies $\
T\models _{GK\leq }\varphi $, hence, $T\models _{GK}\varphi .$
\end{lemma}

Let $T\mathcal{G}_{\Diamond }$ be the set of theorem of $\mathcal{G}%
_{\Diamond }$, then it follows, as in the case of $\mathcal{G}_{\square },$
that

\begin{lemma}
$T\vdash _{\mathcal{G}_{\Diamond }}\varphi $ if and only if $T\cup T\mathcal{%
G}_{\Diamond }\vdash \varphi $ in G\"{o}del logic.
\end{lemma}

Let $\Diamond \mathcal{L}_{\Diamond }=\{\Diamond \theta :\theta \in \mathcal{%
L}_{\Diamond }\}$. The canonical model $\mathcal{M}_{\Diamond }=(W^{\ast
},S^{\ast },e^{\ast })$ is defined as follows:\medskip

\noindent $\bullet $ $\ W^{\ast }$is the set of valuations $v:Var\cup
\Diamond \mathcal{L}_{\Diamond }\rightarrow \lbrack 0,1]$ such that $v(T%
\mathcal{G}_{\Diamond })=1$ and its positive values have a positive lower
bound:%
\begin{equation}
\inf_{\varphi \in \mathcal{L}_{\Diamond }}\{v(\theta ):v(\theta )>0\}=\delta
>0  \label{delta}
\end{equation}%
when the formulas in $\Diamond \mathcal{L}_{\Diamond }$ are seen as
propositional variables and $v$ is extended to $\mathcal{L}_{\Diamond }=%
\mathcal{L}(Var\cup \Diamond \mathcal{L}_{\Diamond })$ as a G\"{o}del
valuation.

\noindent $\bullet $ $\ $The fuzzy relation between worlds in $\mathcal{M}%
_{\Diamond }$ is given by%
\begin{equation*}
S^{\ast }vw:=\inf_{\varphi \in \mathcal{L}_{\Diamond }}\{w(\theta
)\Rightarrow v(\Diamond \theta )\}.
\end{equation*}%
\noindent $\bullet $ $\ e^{\ast }(v,p):=v(p)$ for any $p\in Var$.

\begin{lemma}
{\label{C1}}\emph{For any world }$v$ \emph{in the canonical model }$\mathcal{%
M}_{\Diamond }$\emph{\ and any} $\varphi \in \mathcal{L}_{\Diamond }$ \emph{%
we have }$e^{\ast }(v,\varphi )=v(\varphi ).$
\end{lemma}

\proof The only non trivial step in a proof by induction on complexity of
formulas of $\mathcal{L}_{\Diamond }$ is that of\ $\Diamond .$ By induction
hypothesis, $e^{\ast }(v,\Diamond \varphi )$ $=\sup_{w}\{S^{\ast }vw\cdot
e^{\ast }(w,\varphi )\}$ $=\sup_{w}\{S^{\ast }vw\cdot w(\varphi )\},$ then
we must show $\sup_{w}\{S^{\ast }vw\cdot w(\varphi )\}=v(\Diamond \varphi ).$
By definition 
\begin{equation*}
S^{\ast }vw\leq w(\varphi )\Rightarrow v(\Diamond \varphi ),
\end{equation*}%
for any $\varphi \in \mathcal{L}_{\Diamond }\ $and $w\in W^{\ast },$ then $%
S^{\ast }vw\cdot w(\varphi )\leq v(\Diamond \varphi ),$ which yields taking
join over $w$:%
\begin{equation*}
e^{\ast }(v,\Diamond \varphi )\leq v(\Diamond \varphi ).
\end{equation*}%
The other inequality is trivial if\ $v(\Diamond \varphi )=0.$ For the case $%
v(\Diamond \varphi )>0$, let $w$ be given as in the following claim then $%
v(\Diamond \varphi )=\alpha =S^{\ast }vw\cdot w(\varphi )\leq e^{\ast
}(v,\Diamond \varphi ),$ ending the proof of the lemma.\medskip

\noindent \textbf{Claim}{. }\emph{If }$v$\emph{\ is a world of }$M_{\Diamond
}$\emph{\ such that }$v(\Diamond \varphi )=\alpha >0,$\emph{\ there exists a
world }$w$\emph{\ of }$M_{\Diamond }$\emph{\ such that} $w(\varphi )=1$ and $%
S^{\ast }vw=\alpha .$\medskip

\proof Set 
\begin{equation*}
\Gamma _{\varphi ,v}=\{\theta \in \mathcal{L}_{\Diamond }:v(\Diamond \theta
)<\alpha \}\cup \{\lnot \lnot \mu :\mu \in \mathcal{L}_{\Diamond }\text{, }%
v(\Diamond \mu )=0\}.
\end{equation*}%
This set is not empty because $v(\Diamond 0)=0$ by axiom \textbf{F}$%
_{\Diamond }.$ Moreover, for any finite subset of $\Gamma _{\varphi ,v},$
say $\{\theta _{1},...,\theta _{n}\}\cup \{\lnot \lnot \mu _{1},...,\lnot
\lnot \mu _{m}\},$ we have%
\begin{equation*}
\varphi \nvdash _{_{\mathcal{G}_{\Diamond }}}\theta _{1}\vee ...\vee \theta
_{n}\vee \lnot \lnot \mu _{1}\vee ...\vee \lnot \lnot \mu _{m}.
\end{equation*}%
Otherwise, we would have\medskip

$%
\begin{array}{lll}
\Diamond \varphi \vdash _{\mathcal{G}_{\Diamond }} & \Diamond (\theta
_{1}\vee ...\vee \theta _{n}\vee \lnot \lnot \mu _{1}\vee ...\lnot \lnot
\vee \mu _{m}) & \text{ \textbf{RN}}_{\Diamond } \\ 
& \Diamond \theta _{1}\vee ...\vee \Diamond \theta _{n}\vee \Diamond \lnot
\lnot \mu _{1}\vee ...\vee \Diamond \lnot \lnot \mu _{m} & \text{ \textbf{D}}%
_{\Diamond } \\ 
& \Diamond \theta _{1}\vee ...\vee \Diamond \theta _{n}\vee \lnot \lnot
\Diamond \mu _{1}\vee ...\vee \lnot \lnot \Diamond \mu _{m} & \text{ \textbf{%
Z}}_{\Diamond },%
\end{array}%
$\medskip

\noindent which would imply by Lemma \ref{sound2} 
\begin{equation*}
v(\Diamond \varphi )\leq \max (\{v(\Diamond \theta _{i}):1\leq i\leq n\}\cup
\{v(\lnot \lnot \Diamond \mu _{i}):1\leq i\leq m\})<\alpha ,
\end{equation*}%
a contradiction. Therefore, we have by Lemma \ref{possibilit} 
\begin{equation*}
T\mathcal{G}_{\Diamond },\varphi \nvdash \theta _{1}\vee ...\vee \theta
_{n}\vee \lnot \lnot \mu _{1}\vee ...\vee \lnot \lnot \mu _{m};
\end{equation*}%
\medskip By Proposition \ref{Godel} there is a Heyting algebra valuation $%
u:L\rightarrow \lbrack 0,1]$ such that $u(\varphi )=u(T\mathcal{G}_{\Diamond
})=1$ and $u(\theta )<1$ for all $\theta \in \Gamma _{\varphi ,v}$. Thus,\ $%
u $ satisfies the further conditions:\medskip

\noindent (i) $u(\varphi )=1\newline
$(ii) $u(\theta )<1$\emph{\ if }$v(\Diamond \theta )<\alpha ,$ because then $%
\theta \in \Gamma _{\varphi ,v}\newline
$(iii) $\ u(\theta )=0$\emph{\ if} $v(\Diamond \theta )=0$, because then $%
\lnot \lnot \theta \in \Gamma _{\varphi ,v}$ and so $u(\lnot \lnot \theta
)<1 $ which implies $u(\theta )=0.$ \medskip

\noindent Let $g:[0,1]\rightarrow \lbrack 0,1]$ be the strictly increasing
function:\medskip

$g(x)=\left\{ 
\begin{array}{cc}
1 & \text{if }x=1 \\ 
\delta (x+1)/2 & \text{if }0<x<1 \\ 
0 & \text{if }x=0%
\end{array}%
\right. $\medskip

\noindent where $\delta $ is given by (\ref{delta}). Clearly the valuation $%
w=g\circ u$ inherits the properties (i), (ii)\ (iii) of $u$, with (ii) in
the stronger form:\medskip

(ii$^{\prime }$) $w(\theta )<\delta $ \emph{\ if }$v(\Diamond \theta
)<\alpha $\medskip

\noindent Moreover, $w(\theta )>0$ implies $w(\theta )>\delta /2,$ by
construction, and $w(T\mathcal{G}_{\Diamond })=1$ because $g$ is an
homomorphism of Heyting algebras, hence, $w$ belongs to $\mathcal{M}%
_{\Diamond }.$

To see that $S^{\ast }vw=\alpha ,$ note that $w(\theta )\leq v(\Diamond
\theta )$ whenever $v(\Diamond \theta )<\alpha $. If $0<v(\Diamond \theta )$
because then $w(\theta )<\delta \leq v(\Diamond \theta )$ by (ii$^{\prime }$%
) and definition of $\delta $. If $v(\Diamond \theta )=0$ because then $%
w(\theta )=0$ by (iii). Since $(w(\theta )\Rightarrow v(\Diamond \theta
))\geq \alpha $ for $v(\Diamond \theta )\geq \alpha $, and ($w(\varphi
)\Rightarrow v(\Diamond \varphi ))=(1\Rightarrow \alpha )=\alpha $, we have $%
S^{\ast }vw=\inf_{\varphi \in \mathcal{L}_{\Diamond }}\{w(\varphi
)\Rightarrow v(\Diamond \varphi )\}=\alpha .$\ $\blacksquare $

\begin{theorem}
For any countable theory $T$ and formula $\varphi $ in $\mathcal{L}%
_{\Diamond }$, $T\models _{GK}\varphi $ iff $T\vdash _{\mathcal{G}_{\Diamond
}}\varphi .$
\end{theorem}

\proof Assume that $T\nvdash _{\mathcal{G}_{\Diamond }}\varphi $, then $%
T\cup T\mathcal{G}_{\Diamond }\nvdash \varphi .$\ By strong completeness of G%
\"{o}del logic, there is a Heyting algebra valuation $v$ such that $v(T\cup T%
\mathcal{G}_{\Diamond })=1$ and $v(\varphi )<1.$ Since $v$ might not be a
world in $\mathcal{M}_{\Diamond }$ compose it with the Heyting algebra
homomorphism: $g(x)=(x+1)/2$ for $x>0,$ $g(0)=0.$ Then $v^{\prime }=g\circ v$
belongs to $\mathcal{M}_{\Diamond }$ and we still have $v^{\prime }(T)=1$, $%
v^{\prime }(\varphi )<1.$ Applying Lemma \ref{C1} to $v^{\prime }$ we have $%
e^{\ast }(v^{\prime },T)=1,$ $e^{\ast }(v^{\prime },\varphi )<1.$ That is,%
\emph{\ }$\mathcal{M}_{\Diamond }\models _{v^{\prime }}T$ and $\mathcal{M}%
_{\Diamond }\nvDash _{v^{\prime }}\varphi .$ Hence, $T\nvDash _{\mathcal{G}%
_{\Diamond }}\varphi .$ $\blacksquare $

\medskip

By Lemma \ref{sound} we have again, as in the case of $\mathcal{L}_{\square
},$ that $\models _{GK}$ and $\models _{GK\leq }$ coincide in $\mathcal{L}%
_{\Diamond }.$ However, $\models _{GK}$ no longer coincides with $\models
_{Crisp}$ as the following example illustrates.\medskip

\textbf{Example}. $\mathcal{G}_{\Diamond }$ \emph{is not complete for crisp
models}. The formula $\lnot \lnot \Diamond \varphi \rightarrow \Diamond
\lnot \lnot \varphi $ holds in all crisp models because $e(x,\lnot \lnot
\Diamond \varphi )>0$ implies that there is $y$ such that $Sxy\cdot
e(y,\varphi )>0.$ Thus, $Sxy=1$ and $--e(y,\varphi )=1$ showing that $%
e(y,\Diamond \lnot \lnot \varphi )\geq Sxy\cdot (--e(y,\varphi ))=1$. But
this formula is not a theorem of $\mathcal{G}_{\Diamond }$ because it fails
in the two worlds model: 
\begin{equation*}
x\overset{\frac{1}{2}}{\rightarrow }y,\text{ \ }e(x,p)=e(y,p)=1.
\end{equation*}%
where $e(x,\lnot \lnot \Diamond p)=1,$ and $e(y,\Diamond \lnot \lnot \varphi
)=\frac{1}{2}.$

\section{$\mathcal{G}_{\Diamond }$ has the finite model property}

For any sentence $\varphi $ such that $\nvdash _{_{\mathcal{G}_{\Diamond
}}}\varphi $ we may construct a finite counter-model inside $\mathcal{M}%
_{\Diamond }.$

\begin{theorem}
If $\nvdash _{\mathcal{G}_{\Diamond }}\varphi $\ then there is a model $M$
with finitely many worlds such that $M\nvDash _{GK}\varphi .$
\end{theorem}

\proof It follows from the Claim in Lemma \ref{C1} that for all $\theta $
and $v\in \mathcal{M}_{\Diamond }$ there is $w\in \mathcal{M}_{\Diamond }$
such that $v(\Diamond \theta )=S^{\ast }vw\cdot w(\theta ).$ (if $v(\Diamond
\theta )=0$ any $w$ works). Given $\theta $, let $f_{\theta }(v)$ be a
function choosing one such $w$ for each $v$. For any formula $\theta $ let $%
r(\theta )$ be the \emph{nesting degree} of $\Diamond $\ in $\theta $, that
is, the length of a longest chain of \ occurrences of $\Diamond $ in the
tree of $\theta $.

Given $\varphi $ such that $\nvdash _{\vdash _{\mathcal{G}_{\Diamond
}}}\varphi ,$ let $v_{0}$ be a world (valuation) in $\mathcal{M}_{\Diamond }$
such that $v_{0}(\varphi )<1.$ For each $j\leq n=r(\varphi )$, let $S_{j}$
be the set of subformulas of $\varphi $ of rank $\leq j,$ and define
inductively the following sets of valuations:%
\begin{equation*}
\begin{array}{ccc}
M_{0} & =\{v_{0}\} &  \\ 
M_{i+1} & =M_{i}\cup & \{f_{\theta }(v):v\in M_{i},\text{ \ }\Diamond \theta
\in S_{n-i}\}%
\end{array}%
\end{equation*}%
Clearly, $M_{n}$ is finite. Consider the model induced in $M_{n}$ by
restricting $e^{\ast }$ and $S^{\ast }$ of $\mathcal{M}_{\Diamond }$ to $%
M_{n}\times Var$ and $M_{n}\times M_{n}\ $respectively. We call this model $%
M_{n}$ for simplicity. Then for any formula $\Diamond \theta \in S_{j}$ and $%
v\in M_{n-j}$ there is $w\in M_{n-(j-1)}$ such that $v(\Diamond \theta
)=S^{\ast }vw\cdot w(\varphi ),$ and thus 
\begin{equation*}
v(\Diamond \theta )=\sup_{w}\{S^{\ast }vw\cdot w(\theta ):w\text{ is a world
in }M_{n}\mathcal{\}}.
\end{equation*}%
This permits to show\ by induction in $j\leq n$ that for all $\theta \in
S_{j},$ $v\in M_{n-j}$ we have $v(\theta )=e_{M_{n}}(v,\theta ).$ In
particular, $e_{M_{n}}(v_{0},\varphi )=v_{0}(\varphi )<1,$ which shows $%
M_{n}\nvDash \varphi .$ $\blacksquare $

\bigskip

The proof of the previous theorem still works if we define the accessibility
relation in $M_{n}$ using only subformulas of $\varphi ,$ 
\begin{equation*}
S_{n}^{\ast }vw:=\min_{\theta \in S_{n}}\{w(\theta )\Rightarrow v(\Diamond
\theta )\}
\end{equation*}%
This means that we have to use only a finite number of values of $[0,1]$ in
the proof, and thus $e^{\ast }$ takes values in a finite subalgebra of $%
[0,1].$

\section{Modal extensions}

The modal systems we have considered so far correspond to minimal modal
logic, the logic of G\"{o}del-Kripke models with an arbitrary accessibility
fuzzy relation. We may consider also for each modal operator the analogues
of\ the classical modal systems $T,$ $S4$ and $S5$, usually presented as
combinations of the following axioms:\medskip 
\begin{equation*}
\begin{array}{lllllll}
\mathbf{T}_{\square }\text{:} & \Box \varphi \rightarrow \varphi &  &  & 
\mathbf{T}_{\Diamond }\text{:} & \varphi \rightarrow \Diamond \varphi &  \\ 
4_{\square }\text{:} & \Box \varphi \rightarrow \Box \Box \varphi &  &  & 
4_{\Diamond }\text{:} & \Diamond \Diamond \varphi \rightarrow \Diamond
\varphi &  \\ 
\mathbf{B}_{\square }\text{:} & \varphi \rightarrow \Box \lnot \Box \lnot
\varphi &  &  & \mathbf{B}_{\Diamond }\text{:} & \varphi \rightarrow \lnot
\Diamond \lnot \Diamond \varphi & 
\end{array}%
\end{equation*}%
\medskip

Call a GK-model $\mathcal{M}$ = $\langle W,S,e\rangle $ \emph{reflexive} if $%
Sxx=1$ for all $x\in W$, \emph{(min)transitive} if $Sxy\cdot Syz\leq Sxz$
for all $x,y,z,$ and \emph{symmetric} if $Sxy=Syx$ for all $x,y\in W.$

\begin{proposition}
\label{similarity1}$\mathbf{T}_{\square }$ and $\mathbf{T}_{\Diamond }$ are
valid in all reflexive GK-models, $\mathbf{4}_{\square }$ and $\mathbf{4}%
_{\Diamond }$ are valid in all transitive GK-models, $\mathbf{B}_{\square }$
and $\mathbf{B}_{\Diamond }$ are\textbf{\ }valid in all GK-symmetric models.
\end{proposition}

\proof If $Sxx=1$ for all $x$ then $(\mathbf{T}_{\square })$: $e(x,\Box
\varphi )\leq (Sxx\Rightarrow e(x,\varphi ))=e(x,\varphi ),$ and $(\mathbf{T}%
_{\Diamond })$: $e(x,\Diamond \varphi )$ $\geq Sxx\cdot e(x,\varphi )$ $%
=e(x,\varphi ).$

Assume $Sxy\cdot Syz\leq Sxz$ for all $x,y,z.(\mathbf{4}_{\square })$: $%
e(x,\Box \varphi )\cdot Sxy\cdot Syz\leq (Sxz\Rightarrow e(z,\varphi ))\cdot
Sxz\leq e(z,\varphi ).$ Hence, $e(x,\Box \varphi )\cdot Sxy\leq
(Syz\Rightarrow e(z,\varphi )).$ Taking meet over $z$ in the right hand
side: $e(x,\Box \varphi )\cdot Sxy\leq e(y,\Box \varphi );$ hence, $e(x,\Box
\varphi )\leq (Sxy\Rightarrow e(y,\Box \varphi ))$ for all $y$ and thus $%
e(x,\Box \varphi )\leq e(x,\Box \Box \varphi ).$ ($\mathbf{4}_{\Diamond })$:
For any $x,y,z,$ $Sxy\cdot Syz\cdot e(z,\varphi )$ $\leq Sxz\cdot
e(z,\varphi )\leq e(x,\Diamond \varphi ).$ Hence, $Syz\cdot e(z,\varphi
)\leq (Sxy\Rightarrow e(x,\Diamond \varphi )).$ Taking join over $z$ in the
left$,$ $e(x,\Diamond \varphi )\leq (Sxy\Rightarrow e(x,\Diamond \varphi )),$
thus $Sxy\cdot e(x,\Diamond \varphi )\leq e(x,\Diamond \varphi )).$ Taking
join again in the left, $e(x,\Diamond \Diamond \varphi )\leq e(x,\Diamond
\varphi ).$

Assume $Sxy=Syx$ for all $x,y.$ ($\mathbf{B}_{\square }$): We prove the
stronger $\lnot \varphi \rightarrow \Box \lnot \Box \varphi $. Assume $%
e(x,\lnot \varphi )>0$ then $e(x,\varphi )=0.$ Take any $y$ such that $%
Sxy>0, $ then $e(y,\Box \varphi )\leq (Syx\Rightarrow e(x,\varphi
))=(Sxy\Rightarrow e(x,\varphi ))=0.$ Therefore, $e(y,\lnot \Box \varphi
)=1, $ and ($Sxy\Rightarrow e(y,\lnot \Box \varphi ))=1.$ This shows that $%
x(\Box \lnot \Box \varphi )=1.$ ($\mathbf{B}_{\Diamond })$: Suppose $%
e(x,\varphi )>e(x,\lnot \Diamond \lnot \Diamond \varphi )$ then $e(x,\lnot
\Diamond \lnot \Diamond \varphi )=0$ and $e(x,\Diamond \lnot \Diamond
\varphi )=1.$ This means that there is $y$ such that $Sxy\cdot e(x,\lnot
\Diamond \varphi )>0$ thus $Sxy>0$ and $e(x,\lnot \Diamond \varphi )=1,$
hence $e(y,\Diamond \varphi )=0,$ therefore, $Syx\cdot e(x,\varphi )=0$
which is absurd because $Syx=Sxy>0$ and $e(x,\varphi )>0$ by construction. $%
\blacksquare $\medskip

Let $Ref$, $Trans,$ and $Symm$ denote the GK-classes of models satisfying,
respectively,reflexivity, transitivity, and symmetry, and let $\models _{%
\mathcal{C}}$ denote semantic consequence with respect to models in the
class $\mathcal{C}$.

\begin{theorem}
\label{modales}\emph{(i)}$\ \ \mathcal{G}_{\square }$\emph{+}$\mathbf{T}%
_{\square }$ and $\mathcal{G}_{\Diamond }\emph{+}\mathbf{T}_{\Diamond }$ are
strongly complete (for countable theories) with respect to $\models _{Refl}$.%
\newline
\emph{(ii)} $\mathcal{G}_{\square }\emph{+}\mathbf{4}_{\square }$ and $%
\mathcal{G}_{\Diamond }\emph{+}\mathbf{4}_{\Diamond }$ are strongly complete
with respect to $\models _{Trans}$.\newline
\emph{(iii)} $\mathcal{G}S4_{\square }:=\mathcal{G}_{\square }\emph{+}%
\mathbf{T}_{\square }\emph{+}\mathbf{4}_{\square }$ and $\mathcal{G}%
S4_{\Diamond }:=\mathcal{G}_{\Diamond }\emph{+}\mathbf{T}_{\Diamond }\emph{+}%
\mathbf{4}_{\Diamond }$ are strongly complete with respect to $\models
_{Refl\cap Trans}.$
\end{theorem}

\proof Soundness follows from Proposition \ref{similarity1}. Completeness
follows, in each case, by asking the worlds of the canonical models
introduced in the completeness proofs of $\mathcal{G}_{\square }$ and $%
\mathcal{G}_{\Diamond }$ to satisfy the corresponding schemes. The key fact
is that these schemes force the accessibility relations $S_{\square }^{\ast
}vw=\inf_{\varphi \in \mathcal{L}_{\square }}\{v(\square \varphi
)\Rightarrow w(\varphi )\}$ and $S_{\Diamond }^{\ast }vw=\inf_{\varphi \in 
\mathcal{L}_{\Diamond }}\{v(\varphi )\Rightarrow w(\Diamond \varphi )\}$ to
satisfy the respective properties. (i)\ If $v(\mathbf{T}_{\square })=1$ then 
$S_{\square }^{\ast }vv=\inf_{\varphi \in \mathcal{L}_{\square }}\{v(\square
\varphi \rightarrow \varphi )\}=1.$ If $v(\mathbf{T}_{\Diamond })=1$ then $%
S_{\Diamond }^{\ast }vv=\inf_{\varphi \in \mathcal{L}_{\square }}\{v(\varphi
\rightarrow \Diamond \varphi )\}=1$. (ii)\ If $v(\mathbf{4}_{\square })=1$
then $v(\square \varphi )\leq v(\square \square \varphi )$ and so%
\begin{eqnarray*}
S_{\square }^{\ast }vv^{\prime }\cdot S_{\square }^{\ast }v^{\prime
}v^{\prime \prime } &\leq &[(v(\square \square \varphi )\Rightarrow
v^{\prime }(\square \varphi ))\cdot (v^{\prime }(\square \varphi
)\Rightarrow v^{\prime \prime }(\varphi ))] \\
&\leq &(v(\square \square \varphi )\Rightarrow v^{\prime \prime }(\varphi
))\leq (v(\square \varphi )\Rightarrow v^{\prime \prime }(\varphi ))
\end{eqnarray*}%
Taking meet over $\varphi $ in the last formula we get:\ $S_{\square }^{\ast
}vv^{\prime }\cdot S_{\square }^{\ast }v^{\prime }v^{\prime \prime }\leq
S_{\square }^{\ast }vv^{\prime \prime }.$ (iii)\ If $v(\mathbf{4}_{\Diamond
})=1$ then $v(\Diamond \Diamond \varphi )\leq v(\Diamond \varphi )$ and thus%
\begin{eqnarray*}
S_{\Diamond }^{\ast }vv^{\prime }\cdot S_{\Diamond }^{\ast }v^{\prime
}v^{\prime \prime } &\leq &[(v^{\prime }(\Diamond \varphi )\Rightarrow
v(\Diamond \Diamond \varphi ))\cdot (v^{\prime \prime }(\varphi )\Rightarrow
v^{\prime }(\Diamond \varphi ))] \\
&\leq &(v^{\prime \prime }(\varphi )\Rightarrow v(\Diamond \Diamond \varphi
))\leq (v^{\prime \prime }(\varphi )\Rightarrow v(\Diamond \varphi ))\text{
\ \ }
\end{eqnarray*}%
Minimizing over $\varphi $ in the last formula we get $S_{\Diamond }^{\ast
}vv^{\prime }\cdot S_{\Diamond }^{\ast }v^{\prime }v^{\prime \prime }\leq
S_{\Diamond }^{\ast }vv^{\prime \prime }.$ $\blacksquare $

\medskip

One of the original motivations of the second author to study these fuzzy
modal logics was to interpret the possibility operator $\Diamond $ in the
class of G\"{o}del frames $Refl\cap Trans\cap Symm$ as a notion of
similarity in the sense of Godo and Rodr\'{\i}guez \cite{Godo99}, and a
reasonable conjecture was that $\mathcal{G}S5_{\Diamond }=\mathcal{G}%
S4_{\Diamond }$+$\mathbf{B}_{\Diamond }$ would axiomatize validity in this
frames. Unfortunately, the axioms $\mathbf{B}_{\square }$, $\mathbf{B}%
_{\Diamond }$ do not force symmetry in the canonical models. Thus, we have
not been able to show completeness of $\mathcal{G}_{\square }$+$\mathbf{B}%
_{\square }$ or $\mathcal{G}_{\Diamond }$+$\mathbf{B}_{\Diamond }$ for $%
\models _{Symm}$, even the less completeness of $\mathcal{G}S5_{\Diamond }$
or $\mathcal{G}S5_{\square }=\mathcal{G}S4_{\square }$+$\mathbf{B}_{\square
} $ with respect to $\models _{Refl\cap Trans\cap Symm}$. Perhaps stronger
symmetry axioms as\medskip

$(\varphi \rightarrow \square \theta )\rightarrow \Box (\Box \varphi
\rightarrow \theta )$

$\Diamond (\Diamond \varphi \rightarrow \theta )\rightarrow (\varphi
\rightarrow \Diamond \theta )$\medskip

\noindent would do. In any case, it is possible to show that validity in G%
\"{o}del $Refl\cap Trans\cap Symm$ is decidable.

\section{Adding truth constants}

The previous results on strong completeness may be generalized to languages
with a set \thinspace $Q\subseteq \lbrack 0,1]$ of truth values added as
logical constants to the language, provided $Q$ is topologically discrete
and well-ordered, in particular when $Q$ is finite.

Introduce a constant connective symbol for each $r\in Q,$ denoted by $r$
itself excepting $0$ and $1$ which are identified with $\bot $ and $\top $.
Let $\mathcal{G}_{\square }$($Q)$ be the logic obtained by adding to $%
\mathcal{G}_{\square }$ the axiom schemes R1 - R4 below, and let $\mathcal{G}%
_{\Diamond }(Q)$ be defined similarly by adding to $\mathcal{G}_{\Diamond }$
the book-keeping axioms R1 and R5 - R7, for all $r,s\in Q:$\medskip

R1. (book-keeping axioms)

\ \ \ \ \ $\ r\rightarrow s,$ \ \ \ \ \ if $r\leq s,$

$\ \ \ \ \ \ (r\rightarrow s)\rightarrow s,$ \ \ \ \ if $s<r$

R2. $r\rightarrow \square r$

R3. $(r\rightarrow \square \theta )\rightarrow \square (r\rightarrow \theta
) $

R4. $((\square \theta \rightarrow r)\rightarrow r)\rightarrow \square
((\theta \rightarrow r)\rightarrow r)$\medskip

R5. $\Diamond r\rightarrow r$

R6. $\Diamond (r\rightarrow \varphi )\rightarrow (r\rightarrow \Diamond
\varphi )$

R7. $\Diamond ((\varphi \rightarrow r)\rightarrow r)\rightarrow ((\Diamond
\varphi \rightarrow r)\rightarrow r)$\medskip

\noindent Note that the double negation shift axioms \textbf{Z}$_{\square }$
and \textbf{Z}$_{\Diamond }$ become superfluous (follow from R4 and R7,
respectively). Moreover, R2+R3 may be replaced by the single axiom: $%
(r\rightarrow \square \theta )\leftrightarrow \square (r\rightarrow \theta
), $ and R5+R6 by $\Diamond (r\rightarrow \varphi )\leftrightarrow
(r\rightarrow \Diamond \varphi ).$

The evaluation of GK-models is extended by defining: $e(x,r)=r$ for each $%
r\in Q.$ It may be shown then that $\mathcal{G}_{\square }$($Q)$ and $%
\mathcal{G}_{\Diamond }(Q)$ are strongly complete for countable theories in
their respective languages, and the same holds for the logics mentioned in
Theorem \ref{modales}.

This extends substantially a result of Esteva, Godo and Nogera \cite%
{Esteva04} on weak completeness of G\"{o}del logic with rational truth
constants. If one is interested in weak completeness only, no condition is
needed on $Q$ since it is enough to consider the finitely many truth
constants appearing in the sentence to be proved. On the other hand,
discreteness of $Q$ is necessary for strong completeness:\ if $r$ is a limit
points in $Q$ then there is a strictly increasing or decreasing sequence
converging to $r$, say $\{r_{n}\}$ increases to $\sup r_{n}=r,$ then 
\begin{equation*}
\{r_{1}\rightarrow \theta ,r_{2}\rightarrow \theta ,r_{3}\rightarrow \theta
....\}\models _{GK}r\rightarrow \theta
\end{equation*}%
but no finite subset of premises can grants this, thus by soundness%
\begin{equation*}
\{r_{1}\rightarrow \theta ,r_{2}\rightarrow \theta ,r_{3}\rightarrow \theta
....\}\nvdash _{\mathcal{G}_{\square }(Q)}r\rightarrow \theta ,
\end{equation*}%
Discreteness is not enough, however: $Q=\{r_{1}<r_{2}<....$ \quad $%
....<q_{2}<q_{1}\}$ with $\sup r_{i}=\inf q_{i}$ is discrete, and%
\begin{equation*}
r_{1}\rightarrow \theta ,\text{ }r_{2}\rightarrow \theta ,....\text{ , }\psi
\rightarrow q_{1},\text{ }\psi \rightarrow q_{2},\text{ ....}\models
_{GK}\psi \rightarrow \theta
\end{equation*}%
but no finite subset of the premises yields the same consequence.

\section{Comment}

It rests to axiomatize validity and consequence of the full logic with both
modal operators combined. It may be seen that he union of the systems $%
\mathcal{G}_{\square }$ and $\mathcal{G}_{\Diamond }$ is not enough for that
purpose. However, $\mathcal{G}_{\square }\cup \mathcal{G}_{\Diamond }$
together with Fischer Servi \ \cite{FischerServi84} "connecting axioms":%
\begin{equation*}
\begin{array}{l}
\Diamond (\varphi \rightarrow \psi )\rightarrow (\Box \varphi \rightarrow
\Diamond \psi ) \\ 
(\Diamond \varphi \rightarrow \Box \psi )\rightarrow \Box (\varphi
\rightarrow \psi )%
\end{array}%
\end{equation*}%
may be proved to be a strongly complete axiomatization. This will be shown
in a sequel of this paper.

\end{document}